\documentclass[11pt,reqno]{amsart}

\usepackage[utf8]{inputenc}
\usepackage[margin=1.25in]{geometry}
\usepackage{hyperref}
\usepackage{appendix}
\usepackage{amsfonts}
\usepackage{amssymb}
\usepackage{stmaryrd} 
\usepackage{amsmath}
\usepackage{amsthm}
\usepackage{mathrsfs}

\theoremstyle{plain}
\newtheorem{theorem}{Theorem}[section]
\newtheorem{corollary}[theorem]{Corollary}

\newtheorem{proposition}[theorem]{Proposition}
\theoremstyle{remark}

\newtheorem{remark}[theorem]{Remark}

\usepackage{biblatex} 
\numberwithin{equation}{section}
\title[Support SPDE Hölder coefficients]{A support theorem for parabolic stochastic PDEs with nondegenerate Hölder diffusion coefficients}

\author{Yi HAN}
\address{Department of Pure Mathematics and Mathematical Statistics, University of Cambridge.
}
\email{yh482@cam.ac.uk}
\thanks{Supported by EPSRC grant EP/W524141/1.}

\begin{document}

\maketitle

\begin{abstract}
    In this paper we work with parabolic SPDEs of the form
    $$ \partial_t u(t,x)=\partial_x^2 u(t,x)+g(t,x,u)+\sigma(t,x,u)\dot{W}(t,x) 
    $$ with Neumann boundary conditions, where $x\in[0,1]$, $\dot{W}(t,x)$ is the space-time white noise on $(t,x)\in[0,\infty)\times [0,1]$, $g$ is uniformly bounded, and the solution $u\in\mathbb{R}$ is real valued. The diffusion coefficient $\sigma$ is assumed to be uniformly elliptic but only Hölder continuous in $u$. Previously, support theorems for SPDEs have only been established assuming that $\sigma$ is Lipschitz continuous in $u$. We obtain new support theorems and small ball probabilities in this $\sigma$ Hölder continuous case via the recently established sharp two sided estimates of stochastic integrals.
\end{abstract}

\section{Introduction}

The support theorem for stochastic processes has a long history. One of its simplest forms can be phrased as follows: let $B_t$ be a $d$-dimensional Brownian motion started from $0$, then for any $\epsilon,t>0$ we have 
$$\mathbb{P}(\sup_{s\leq t}|B_s|<\epsilon)>0.$$ This follows from the reflection principle of Brownian motion. Via a Girsanov change of measure, we can deduce that for any continuous $\psi:[0,t]\to\mathbb{R}^d$ with $\psi(0)=0$, we have 
$$\mathbb{P}(\sup_{s\leq t}|B_s-\psi(s)|<\epsilon)>0.$$ See for example \cite{bass1994probabilistic}, page 59-60, (6.5) and (6.6). Both claims rely on Gaussian structure of the process $B$.

There is another notion, usually named Stroock–Varadhan support theorem for stochastic processes, that has a different flavour. Consider the parabolic SPDE with periodic boundary conditions
\begin{equation}
    \label{eq*}
\partial_t u(t,x)=\partial_x^2 u(t,x)+g(t,x,u)+\sigma(t,x,u)\dot{W}(t,x),\quad x\in[0,1].
\end{equation} Denote by $\mathcal{H}:=\{h:[0,T]\times[0,1]\to\mathbb{R}, h \text{ absolutely continuous }, \dot{h}\in L^2([0,T]\times[0,1])\}$, the Cameron-Martin space of the Brownian sheet. $\mathcal{H}$ is a Hilbert space endowed with the norm $\|h\|_{\mathcal{H}}:=\|\dot{h}\|_{L^2([0,T]\times[0,1])}$. Let $S(h)$ denote the solution of \eqref{eq*} when we take  $\dot{h}$ in place of the white noise $\dot{W}$ for the given $h\in\mathcal{H}$.
If we assume $\sigma$ is Lipschitz, $u(0,\cdot)$ is Hölder continuous and $g$ has sufficiently many derivatives, it is proved in \cite{bally1995approximation} that the topological support of the probability law of the random variable $u$ in $\mathcal{C}^\alpha([0,T]\times[0,1])$ (for some $\alpha\in(0,\frac{1}{4})$ ) is given by the closure of $\cup_{h\in\mathcal{H}}S(h)$. Here $\mathcal{C}^\alpha([0,T]\times[0,1])$ is the Wiener space equipped with the $\alpha$-Hölder topology.

Support theorems of this flavour have been proved for other examples of SPDEs, starting from \cite{millet1994support}, and see also \cite{cardon2001support}, \cite{chueshov2011stochastic}  and \cite{delgado2014approximation}, with the same conclusion that the topological support of the solution $u$ is the closure of $S(h)$, the solution of the SPDE driven by the control $h$ in place of the noise. The support theorem for singular SPDEs has been obtained in \cite{hairer2022support}, with the feature that the coset structure in the renormalization group plays a key role in characterizing topological supports of the solution. In all these works, the coefficients $g$ and $\sigma$ are nice enough so that the control problem $S(h)$ can be properly solved, and the support of \eqref{eq*} can be characterized in terms of solutions to the control problem $S(h)$.

In this paper we consider support theorems of SPDEs in the lens of the regularization by noise phenomenon. In the setting of \eqref{eq*}, this means that the coefficients $g$ or $\sigma$ (or both) are not necessarily locally Lipschitz continuous. One usually requires $\sigma$ to be uniformly elliptic, so that roughness of the driving noise can restore well-posedness of the equation. We generally do not have a Stroock–Varadhan support theorem since the ODE or PDE for the control process $S(h)$ is in general not well-posed. However it is still possible to prove the solution has full support or obtain small ball probability estimates

For additive noise, i.e., $\sigma=I_d$, and assuming the drift $g$ is not too singular, we can remove the drift via Girsanov transform and show the solution has full support because white noise has. A related example for finite dimensional SDEs with singular drifts can be found in \cite{ling2022wong}. The story is the same when $\sigma$ is Lipschitz continuous and uniformly elliptic.

We are particularly interested in the \textit{robustness} of the support theorems, in the remaining (hardest) case that $\sigma$ is only $\alpha$-Hölder continuous in $u$, for $\alpha\in(0,1]$. More precisely, we assume that for some $\mathcal{C}_1,\mathcal{C}_2,\mathcal{D}>0$ we have 
$$|\sigma(t,x,u)-\sigma(t,x,v)|\leq\mathcal{D}|u-v|^\alpha,\quad\mathcal{C}_1\leq |\sigma(t,x,u)|\leq\mathcal{C}_2.$$ for any $t>0$, $x\in[0,1]$ and $u,v\in\mathbb{R}$.
In this case the Girsanov transform or the control processes $S(h)$ do not tell us the answer in the same way as before.

Before stating our support theorems, it is crucial to discuss the well-posedness issues of \eqref{eq*}. For $\sigma$ to be $\alpha$-Hölder continuous in $u$, $\alpha>\frac{3}{4}$, strong existence and strong uniqueness have been established in \cite{mytnik2011pathwise} (see also \cite{han2022exponential} for a different perspective, also with $\alpha>\frac{3}{4}$). In general we may consider probabilistic weak solutions to \eqref{eq*} without the knowledge of weak uniqueness (weak existence follows from \cite{gkatarek1994weak}). So long as $\alpha>0$, the lower bound in our main theorem still holds, and in particular every weak solution to \eqref{eq*}, considered as a probability law on $\mathcal{C}([0,T]\times[0,1])$ (with the supremum norm),  has full topological support in the sense that the support is given by $\{w(t,x)\in\mathcal{C}([0,T]\times[0,1]):w(0,x)=u(0,x)\}$.

Our strategy to characterize the support of \eqref{eq*} goes as follows. Since $\sigma$ is non-degenerate, one expects that \eqref{eq*} lies in the same universality class as the linear stochastic heat equation 
\begin{equation}\label{linear12}\partial_t u(t,x)=\partial_x^2u(t,x)+\dot{W}(t,x),\end{equation}
which is also named as the Edwards-Wilkinson universality class. In this universality class we observe nontrivial limiting behavior under the 1:2:4 scaling relation $u(t,x)\mapsto \epsilon^{-1}u(\epsilon^{-4}t,\epsilon^{-2}x)$.\footnote{\label{footnote1} 
Strictly speaking, the 1:2:4 scaling does not readily apply to SPDEs on the spatial domain $[0,1]$, but we will use the following more flexible interpretation of 1:2:4 scaling: consider the stochastic integral with respect to the Brownian sheet on small scales, more precisely on a spatial scale $[0,\epsilon^2]$ and temporal scale $[0,\epsilon^4]$, then the stochastic integral typically has values on a scale $[-\epsilon,\epsilon]$ with good tail estimates (See \cite{athreya2021small}, Lemma 3.4). We then derive from these small-scale estimates a macroscopic estimate.} From this scaling relation we can expect fairly sharp two-sided probability estimates for the solution of \eqref{eq*} on small scales, even in the case that $\sigma$ is not constant in $u$. Such computations have been carried out in \cite{athreya2021small}, obtaining matching small ball probability estimates of solutions to \eqref{eq*} when $\sigma$ is Lipschitz continuous in $u$, with a Lipschitz constant sufficiently small.

When $\sigma$ is merely Hölder continuous, we expect that it induces a highly nonlinear stretching of space and time, in the following sense: if $u$ is approximated by another Gaussian field $u'$ on a space-time grid that preserves the 1:2 parabolic scaling, when $\sigma$ is Lipschitz with small Lipschitz constant we can prove that $u$ stays close to $u'$ on the desired scale via stochastic calculus, but when $\sigma$ is merely Hölder we cannot prove $u$ lives close to $u'$ on that scale with high probability. To overcome this, we adjust the 1:2 space-time parabolic scaling through a reduction of the temporal length scale while keeping the spatial scale fixed. Consequently, we can still obtain nontrivial (upper and lower) probability estimates of fine scale properties of the solution, and this is already sufficient for us to prove the support theorem. The upper and lower bounds of small ball probabilities in this Hölder continuous case however do not have matching exponents in $\epsilon$, reflecting the fact that irregularity of $\sigma$ induces high order stretching in space and time. We note that for the linear stochastic heat equation \eqref{linear12}, we can obtain small ball probabilities where the upper and lower bounds have matching exponents in $\epsilon,$ see \cite{dalang2009minicourse}, page 168, Theorem 5.1.

We now state the main theorem.

Denote by $\mathcal{P}$ the predictable $\sigma$-field of the noise $\dot{W}(t,x)$, generated by functions of the form $f(x,t,\omega)=X(\omega)\cdot 1_A(x)\cdot 1_{(a,b]}(t),$ with $A\subset [0,1]$ and $X$ some $\mathcal{F}_a$-measurable random variable. We say $h\in\mathcal{P}\mathcal{C}_b^2$ if $h\in\mathcal{P}$ and almost surely, $h,\partial_t h,\partial_x^2h$ are bounded by a fixed constant $C$.

An important remark before the statement: assuming $\sigma$ is $\alpha$-Hölder in $u$ for any $\alpha>0$, then by a compactness argument there always exists probabilistic weak solutions to \eqref{asterisque}, see for example \cite{gkatarek1994weak}. However, there is no proof in the literature that the solution is unique for general $\alpha$. See Remark \ref{remark1.3} for pathwise well-posedness results for $\alpha$ in certain regimes.

\begin{theorem}\label{theorem1.1}
Consider solution $u(t,x)\in\mathbb{R}$ to the stochastic heat equation with Neumann boundary conditions on $[0,1]$
\begin{equation}\label{asterisque}
\partial_t u(t,x)=\frac{1}{2}\partial_x^2 u(t,x)+g(t,x,u(t,x))+\sigma(t,x,u(t,x))\dot{W}(t,x),\quad u(0,x)=u_0(x).\end{equation}
Assume $u_0,h\in\mathcal{P}\mathcal{C}_b^2$ with $$\sup_{x\in[0,1]}|u_0(x)-h(0,x)|<\epsilon/2,$$ and that for some constants $D,\mathcal{C}_1,\mathcal{C}_2>0$, $\alpha\in(\frac{3}{4},1],$ we have
$$|\sigma(t,x,u)-\sigma(t,x,v)|\leq \mathcal{D}|u-v|^\alpha,$$
$$\mathcal{C}_1\leq |\sigma(t,x,u)|\leq \mathcal{C}_2$$ for all $x\in[0,1]$, $u,v\in\mathbb{R}$ and $t\geq 0$, and 
$$
\sup_{t>0,x\in[0,1],u\in\mathbb{R}}|g(t,x,u)|<\infty
.$$
Then for any $\beta>2-\alpha$ we may find positive constants $C_0,C_1,C_2,C_3$ and $\epsilon_0$ depending on $\beta,\mathcal{C}_1,\mathcal{C}_2$ and $\sup_{t,x,u}|g(t,x,u)|$, such that for any $0<\epsilon<\epsilon_0$, we have
\begin{equation}\label{mainest}
C_0\exp(-\frac{C_1 T}{\epsilon^{2+4\beta}})\leq P(\sup_{0\leq t\leq T,x\in[0,1]}|u(t,x)-h(t,x)|\leq\epsilon)\leq C_2\exp(-\frac{C_3 T}{(1+\mathcal{D}^2)\epsilon^{4+2\alpha}}).\end{equation}
If we only assume $\alpha>0$, then the lower bound in \eqref{mainest} holds, that is,
\begin{equation}\label{1.sgas}
    C_0\exp(-\frac{C_1 T}{\epsilon^{2+4\beta}})\leq P(\sup_{0\leq t\leq T,x\in[0,1]}|u(t,x)-h(t,x)|\leq\epsilon).\end{equation}
Moreover, there exists a $\mathcal{D}_0>0$ depending only on $\mathcal{C}_1$, $\mathcal{C}_2$ such that whenever $\mathcal{D}<\mathcal{D}_0$, the above estimates \eqref{mainest} and \eqref{1.sgas} can hold for $\beta=2-\alpha$ (the value of the various numerical constants may be changed.)
\end{theorem}
Since $h$ is arbitrary, this in particular implies that the solution $u$ has full support on Wiener space with respect to the supremum norm. More precisely:

\begin{corollary} Assume $\alpha\in(0,1]$.
Let $\mu$ denote the probability law of any possible weak solution to the SPDE \eqref{asterisque} on $\mathcal{C}([0,T]\times[0,1]).$ Denote by 
$$  \Omega_{u_0}([0,T]\times[0,1])=\left\{w(t,x)\in \mathcal{C}([0,T]\times[0,1]): w(0,x)=u_0(x)\right\},$$
then for any such $\mu$ there holds
$$\operatorname{Supp}\mu=\Omega_{u_0}([0,T]\times[0,1]),$$ 
where $\mathcal{C}([0,T]\times[0,1])$ is endowed with the topology generated by the supremum norm.
\end{corollary}

Theorem \ref{theorem1.1} is not yet satisfactory in that the upper and lower bounds in \eqref{mainest} have a very wide gap and is thus likely to be sub-optimal. This is to be compared with the main result in \cite{athreya2021small}
,where the upper and lower estimates have matching powers of $\epsilon$ whenever the Lipschitz constant of $\sigma$ is sufficiently small. We believe that the potentially sub-optimal estimate \eqref{mainest} arises from purely technical limitations, but these technical constraints are hard to remove in our infinite dimensional setting. Indeed, when one deals with finite dimensional diffusion with no drift, one can represent the solution as a time changed Brownian motion, and the support theorem readily follows. This is however not the case for SPDEs, and for SPDEs we usually can only approximate the solution by some Gaussian variables, as done in \cite{athreya2021small} and this paper. (This also possibly explains why small ball probability estimates for finite dimensional diffusion are long known and easy to prove but the corresponding estimates for SPDEs are only obtained recently.)  As discussed in footnote \ref{footnote1}, on small scales a 1d Gaussian field is expected to have a 1:2:4 scaling, but if we compose this Gaussian field by a Hölder continuous function  which is not Lipschitz, the 1:2:4 scaling will be distorted. Our main interest in this fine scale property is that we will approximate the solution $u$ by some Gaussian field $u'$, and if this optimal 1:2:4 scaling is violated, our proof only gives an upper bound of distance of $u$ and $u'$ that is much larger in magnitude than what is considered to be optimal. To sum up, this microscopic, approximation by Gaussian procedure is precisely the place that needs the 1:2:4 scaling and thus will lead to suboptimal estimates for Hölder continuous $\sigma$. If we can prove the support theorem and small ball probability estimate without using this approximation by Gaussian procedure, we might obtain a sharper, or even matching estimate that improves \eqref{mainest}.

\begin{remark}\label{remark1.3}
The technical assumption $\alpha\in(\frac{3}{4},1]$ is only used to match with the well-posedness results in \cite{mytnik2011pathwise} (see also \cite{han2022exponential}). 

 For the upper bound to be proved in Section \ref{section2.4}, we are not sure if the proof carries over for any $\alpha>0$ or not because we have to solve another SPDE \eqref{auxiliary1} \footnote{with the diffusion coefficient $\sigma$ only $\alpha$- Hölder continuous in the solution argument} on the same probability space. This procedure requires strong well-posedness for every diffusion coefficient with the given Hölder continuity.

\end{remark}
\begin{remark}
We have for simplicity worked on the unit interval $[0,1]$, but everything carries over to finite intervals $[0,J]$ for any $J>0$. In this paper we assume the solution $u$ is real-valued, but at least the lower bound \eqref{1.sgas} carries over to the vector valued case $\textbf{u}(t,x)\in\mathbb{R}^d$, $d\in\mathbb{N}_+$ without change. These extensions can be found in \cite{athreya2021small} when $\sigma$ is Lipschitz. The upper bound in \eqref{mainest} might be hard to extend to higher dimensions as the corresponding pathwise well-posedness results are lacking.
\end{remark}

There are a few remaining questions. One of them is to obtain support theorems in Hölder semi-norm rather than the supremum norm. Fairly sharp results have been obtained when $\sigma$ is Lipschitz continuous in $u$ (see \cite{foondun2022small} for recent progress), but adapting estimates in the existing literature to our Hölder continuous $\sigma$ seems a bit out of reach. Another more fundamental question is if we assume $\sigma$ is uniformly elliptic, it is not clear whether the assumption that the Hölder index $\alpha>\frac{3}{4}$ is necessary for (strong or weak )well-posedness of \eqref{asterisque}. We believe that $\alpha>0$ is enough for weak well-posedness but could not give a proof. Note that when $\sigma$ is not assumed to be uniformly elliptic, then $\alpha>\frac{3}{4}$ is a sharp condition, see \cite{mueller2014nonuniqueness} for the case $\alpha<\frac{3}{4}.$

\section{Proof of main theorem}

\subsection{Reduction to simple cases} We quote relevant reduction steps from \cite{athreya2021small}, Section 2.2 for sake of completeness.
We first show that after some simple reductions we can assume $g=0$ and $h=0$. These reductions follow from Girsanov theorem and the non-degeneracy of $\sigma$.

The solution to \eqref{asterisque} can be reformulated as 
$$\partial_t u(t,x)=\frac{1}{2}\partial_x^2 u(t,x)+\sigma(t,x,u(t,x))[\dot{W}(t,x)+\sigma^{-1}g(t,x,u(t,x))].$$

For each $t>0$ let $P_t$ be $\mathcal{P}$ restricted to the filtration $\mathcal{F}_t$.
Consider the probability law $Q_t$ defined as 
\begin{equation}\label{fb}\frac{dQ_t}{dP_t}=\int_0^t \int_0^1 \sigma^{-1}g(t,x,u(t,x)) \cdot W(t,x)-\frac{1}{2}\int_0^t\int_0^1 (\sigma^{-1}g(t,x,u(t,x))^2 dxdt.\end{equation}

By Girsanov theorem, $\dot{W}(t,x)+\sigma^{-1}g(t,x,u(t,x))$ is a space-time white noise with respect to $Q_t$. Denote by $A$ the event that 
$$A=\{\sup_{s\in[0,T],y\in[0,1]}|u(s,y)-h(s,y)|<\epsilon\},$$ then by Cauchy-Schwartz inequality
\begin{equation}\label{cs}Q_T(A)=E^{p_T}[1_A\frac{dQ_T}{dP_T}]\leq \sqrt{P_T(A)}\sqrt{E(\frac{dQ_T}{dP_T})^2}\leq\sqrt{P_T(A))}M,\end{equation}
where $M$ depends only on $T,\mathcal{C}_1$ and $\sup_{t,x,y}|g(t,x,y)|.$ This implies the lower bound in \eqref{mainest} for general $g$ can be deduced from the lower bound in the case $g=0$.

For the upper bound, a similar argument holds: one only needs to swap $P$ and $Q$ in \eqref{cs}
and replace $g$ by $-g$ in \eqref{fb}.

Now we show why we can take $h=0$ . This is outlined in page 6 of \cite{athreya2021small} but we reproduce here for completeness.
Let $H:=\partial_t-\frac{1}{2}\partial_x^2$ and consider the process $$w(t,x)=u(t,x)-u_0(x)-h(t,x)+h_0(x),$$ so that $w(0,x)=0.$ If we set $\sigma_1(t,x,w)=\sigma(t,x,u)$ and $$g_1(t,x,w)=g(t,x,u)-Hu_0(x)-Hh(t,x)+Hh_0(x),$$ we have 
$$\partial_t w(t,x)=\frac{1}{2}\partial_x^2 w(t,x)+g_1(t,x,w)+\sigma_1(t,x,w)\cdot \dot{W}(t,x).$$
Since $u_0,h\in\mathcal{P}\mathcal{C}_b^2,$ $\sup_{t,x,\omega}|g_1(t,x,\omega)|<\infty$, so we are reduced to the case $h=0$ and $u(0,x)\equiv 0,x\in[0,1]$.

\subsection{Sharp two-sided estimates}

Recall the heat kernel on $[0,1]$ is given by 
$$G(t,x)=\sum_{n\in\mathbb{Z}}(2\pi t)^{-1/2}\exp(-\frac{(x+n)^2}{2t}).$$
Consider the noise term $\mathbf{N}$ defined as
$$\mathbf{N}(t,x):=\int_0^t\int_0^1 G(t-s,x-y)\sigma(s,y,u(s,y))W(dyds).$$

We quote the following large deviations estimate from \cite{athreya2021small}, Proposition 3.4 and Remark 3.1, which is a very precise formulation of the 1:2:4 scaling of Gaussian processes:

\begin{proposition}\label{prop2.1}
Assume that $\sup_{s,y}|\sigma(s,y,u(s,y))|\leq\mathcal{C}<\infty$. Then we can find universal constants $K_1$ and $K_2$ such that, for any $\alpha,\lambda,\epsilon>0$, 

\begin{equation}\label{estimate2.11}
\mathbb{P}\left(\sup_{0\leq t\leq \alpha\epsilon^4,x\in[0,\epsilon^2]}|\mathbf{N}(t,x)|>\lambda\epsilon\right)\leq\frac{K_1}{1\wedge \sqrt{\alpha}}\exp(-K_2\frac{\lambda^2}{\mathcal{C}^2\sqrt{\alpha}}).
\end{equation}
\end{proposition}

\begin{remark}\label{remark2.2}
    The estimate \eqref{estimate2.11} also holds if the supremum is taken over $x\in[(k-1)\epsilon^2,k\epsilon^2]$ for all $k\in\mathbb{N}_+$ such that $k\epsilon^2\leq 1$. This is because the proof of \eqref{estimate2.11} in \cite{athreya2021small} only uses a modulus of continuity estimate of $\mathbf{N}(t,x)$ in $t$ and $x$ (check for \cite{athreya2021small}, Lemma 3.3), and this estimate is independent of the location of $x$ in $[0,1]$. Also used is a path of grid points from $(0,0)$ to $(t,x)$, which in this case can be replaced by a path starting from $(0,(k-1)\epsilon^2)$.
\end{remark}

We fix a sufficiently small $c_0>0$ such that $0<c_0<\max\{(\frac{K_2}{36\log K_1 \mathcal{C}_2^2})^2,1\}$, and define the discretized mesh of time as: \footnote{Later we will introduce a different scheme to divide time intervals, which doesn't follow the 1:2 parabolic scaling.} 
 $$t_n=nc_0\epsilon^4,\quad n\geq 0,$$ 
and denote by $I_n:=[t_n,t_{n+1}]$ the time interval with numbering $n$. Choose some $\theta=\theta(\mathcal{C}_1,\mathcal{C}_2)>0$ sufficiently large (with the precise condition given in \cite{athreya2021small},  (2.11)) and fix $c_1^2=\theta c_0$, the spatial mesh points  are chosen as 
$$x_n=nc_1\epsilon^2,n\geq 0.$$

This time-space mesh respects the parabolic 1:2 scaling. Fix a terminal time $T>0$ and define the terminal index 
$$n_1:=\min\{n\geq 1:t_n>T\},\quad n_2:=\min\{n\geq 1:x_n>1\}.$$

Write $p_{i,j}=(t_i,x_j)$, and consider the following two series of events
$$A_n=\{|u(t_{n+1},x)|\leq\frac{\epsilon}{3}, \quad x\in[0,1], \text{ and } |u(t,x)|\leq\epsilon ,t\in I_n,x\in[0,1]\},
$$
and 
$$F_n=\{|u(p_{nj})|<\epsilon\text{ for all }j\leq n_2-2\}.$$

The strategy of proof is first to fix the $u$ component of $\sigma$ and obtain an estimate in the Gaussian case, then deduce the general case via an interpolation argument. For the Gaussian case (when $\sigma$ does not depend on $u$), we quote the following result from \cite{athreya2021small}, Proposition 2.1:

\begin{proposition}\label{prop2.2}
Under the assumptions of Theorem \ref{theorem1.1}, assume further that $g=0$, $u_0(x)\equiv0$ and $\sigma(t,x,u)$ does not depend on $u$. 

Then there exists constants $\epsilon_0,C_4,C_5>0$ which depend only on $\mathcal{C}_1$ and $\mathcal{C}_2$ such that for any $0<\epsilon<\epsilon_0,$
$$P(F_n\mid\cap_{k=0}^{n-1}F_k)\leq C_4\exp(-C_5\epsilon^{-2}),$$
 and we can find constants $C_6,C_7>0$ which depend only on $\mathcal{C}_1,\mathcal{C}_2$ such that for any $0<\epsilon<\epsilon_0,$
    $$P(A_n\mid\cap_{k=0}^{n-1}A_k)\geq C_6\exp(-C_7\epsilon^{-2}).$$
\end{proposition}

Now we prove the general case (i.e., $\sigma$ depends on $u$).
\subsection{Upper bound, general case}\label{section2.4}
Define a function 
$$f_\epsilon(z)=\begin{cases}z,\quad |z|<\epsilon,\\\frac{\epsilon}{|z|}z,\quad |z|>\epsilon,\end{cases}$$
so that $|f_\epsilon(z)|\leq\epsilon$ and $f_\epsilon$ is Lipschitz continuous. We solve the following SPDE
\begin{equation}\label{auxiliary1}\partial_t v(t,x)=\frac{1}{2}\partial_x^2 v(t,x)+\sigma(t,x,f_\epsilon(v(t,x)))\cdot \dot{W}(t,x)\end{equation} with $v(0,x)=u_0(x)$, which is well posed because $\sigma(t,x,f_\epsilon(u))$ is $\alpha$-Hölder continuous in $u$, for $\alpha>\frac{3}{4}$. 

As long as $|u(t,x)|\leq\epsilon$ for all $x\in[0,1]$ and $t\in[0,t_1]$, we have $v(t,x)=u(t,x)$, so we proceed with the proof for $v$. 

The point is to compare $v$ with an auxiliary process $v_g$ defined by
$$
\partial_t v_g(t,x)=\frac{1}{2}\partial_x^2 v_g(t,x)+\sigma(t,x,f_\epsilon(u_0(x)))\cdot \dot{W}(t,x),
$$ with $v_g(0,x)=u_0(x)$,
where the diffusion coefficient is independent of $v_g$. The subscript $g$ in the notation $v_g$ stands for Gaussian.

The difference process $D(t,x):=v(t,x)-v_g(t,x)$ is a stochastic integral satisfying 
$$
D(t,x)=\int_0^t\int_0^1 G(t-s,x-y)[\sigma(s,y,f_\epsilon(v(s,y))-\sigma(s,y,f_\epsilon(u_0(y))]\cdot W(dyds).$$
Define 
$$H_j=\{v(p_{1j})|\leq\epsilon\}$$
and consider the events 
$$A_{1,j}=\{|v_g(p_{1j})|\leq 2\epsilon\},$$
$$A_{2,j}=\{|D(p_{1j})|\geq\epsilon\}.$$
It is clear that $H_j\subset A_{1,j}\cup A_{2,j}$.
Define another sequence of events 
$$B_n=\{|u(t,x)|\leq\epsilon,t\in I_{n-1},x\in [0,1].\},\quad n\geq 1.$$ Then clearly $B_n\subset F_n$ and on $B_n$, we have $u(t_n,x)=v(t_n,x)$. Therefore 
$$P(B_1)\leq P(\cap_{j=1}^{n_2-2}H_j)\leq P(\cap_{j=1}^{n_2-2}(A_{1,j}\cup A_{2,j})).$$
An elementary set-inclusion argument implies
$$P(B_1)\leq P(\cap_{j=1}^{n_2-2}A_{1,j})+\sum_{j=1}^{n_2-2}P(A_{2,j}).$$
We now apply Proposition \ref{prop2.2} to the process $v_g$ to deduce 
$$P(\cap_{j=1}^{n_2-2}A_{1,j})=P(|v_g(p_{1,j})|\leq2\epsilon, \quad j=1,\cdots,n_2-2)\leq C_2\exp(-C_3\epsilon^{-2}).$$

By Hölder continuity of $\sigma$ in $u$, we deduce that 
$$|\sigma(s,y,f_\epsilon(v(s,y)))-\sigma(s,y,f_\epsilon(v_0(y)))|\leq \mathcal{D}(2\epsilon)^\alpha,$$ so that by Proposition \ref{prop2.1} and Remark \ref{remark2.2}, we have for $j=1,\cdots,n_2-2,$
$$P(A_{2,j})\leq K_1\exp(-\frac{K_2}{4\epsilon^{2\alpha}\mathcal{D}^2\sqrt{c_0}}).$$
Therefore 
$$\begin{aligned}
P(B_1)&\leq  C_2\exp(-C_3\epsilon^{-2})+\sum_{j=1}^{n_2-2}K_1\exp(-\frac{K_2}{4\epsilon^{2\alpha}\mathcal{D}^2\sqrt{c_0}})\\& \leq  C_2\exp(-C_3\epsilon^{-2})+\frac{1}{c_1\epsilon^2}K_1\exp(-\frac{K_2}{4\epsilon^{2\alpha}\mathcal{D}^2\sqrt{c_0}})\\&\leq C_4\exp\left(-\frac{C_5}{8(1+\mathcal{D}^2)\epsilon^{2\alpha}\sqrt{c_0}}\right),
\end{aligned}$$
whenever $\epsilon$ is small enough, the $\epsilon^{-2\alpha}$ term in the exponent wins over the $\epsilon^{-2}$ term, so we keep the former. \footnote{We have used another approximation which follows from the elementary inequality $$\sup_{x\geq 0}x^n e^{-Cx}<\infty$$ for any $C>0$ and $n>0$.} $C_4$ and $C_5$ are universal constants that depend only on $\mathcal{C}_2$.

The expression shows that when $\sigma$ is merely Hölder continuous, i.e. $\alpha<1$, the $\epsilon^{2\alpha}$ term dominates in the upper bound. The upper and lower bounds we obtain will not have matching exponents of $\epsilon$ (they do if $\alpha=1$), but both bounds are nontrivial and in particular they lead to the desired support theorem.

By the Markov property, for each $n\leq n_1,$ $$P(B_n\mid\cap_{j=1}^{n-1}B_j)\leq \exp(-\frac{C_7}{(1+\mathcal{D}^2)\epsilon^{2\alpha}})$$
where $C_7$ depends only on $\mathcal{C}_1,\mathcal{C}_2$ and $\mathcal{D}$.

Therefore 
$$\begin{aligned}
P(|u(t,x)|\leq\epsilon,t\in[0,T],x\in[0,1])&=P(\cap_{j=1}^{n_1-1}B_j)\\&\leq \exp(-\frac{\mathcal{C}_7}{(1+\mathcal{D}^2)\epsilon^{2\alpha}})^{\frac{T}{\epsilon^4}}
\\&\leq \exp(-\frac{C_7 T}{(1+\mathcal{D}^2)\epsilon^{2\alpha+4}}). 
\end{aligned}$$
This establishes the upper bound in \eqref{mainest}.
\subsection{Lower bound, general case}\label{section2.5}
We now proceed to prove the corresponding lower bound. The argument roughly follows that in \cite{athreya2021small}, while the last key estimates are different. 

Fix some $\beta>1$ to be determined later, consider a new time mesh as follows:
$$\hat{t}_n:=nc_0\epsilon^{4\beta},n\geq 0,$$ and the corresponding time intervals $\hat{I}_n:=[\hat{t}_n,\hat{t}_{n+1}]$. This introduces a finer grid of time when $\epsilon$ is sufficiently small. We analogously define the events 
$$\hat{A}_n:=\{|u(\hat{t}_{n+1},x)|\leq\frac{\epsilon}{3}, \quad x\in[0,1], \text{ and } |u(t,x)|\leq\epsilon ,t\in \hat{I}_n,x\in[0,1]\},
$$
Assuming that $|u_0(x)|\leq\frac{\epsilon}{3}$, $x\in[0,1]$. Define the stopping time $$\tau=\inf\{t\geq 0:\sup_{x\in[0,1]}|u(t,x)-u_0(x)|>2\epsilon\},$$ such that on the event $\hat{A}_0$ we have $\tau\geq \hat{t}_1.$ Consider the process$$
\widetilde{D}(t,x)=\int_0^t\int_0^1 G(t-s,x-y)[\sigma(s,y,u(s\wedge\tau,y))-\sigma(s,y,u_0(y))]W(dyds),$$
and the auxiliary comparison process $u_g$ solving
$$ \partial_t u_g=\frac{1}{2}\partial_x^2 u_g+\sigma(t,x,u_0(x))dW(t,x),\quad u_g(0,x)=u_0(x),
$$
and we write $u(t,x)=u_g(t,x)+D(t,x)$ as before, so that 
$$
D(t,x)=\int_0^t\int_0^1 G(t-s,x-y)[\sigma(s,y,u(s,y))-\sigma(s,y,u_0(y))]W(dyds).$$

It is clear that $D(t,x)=\widetilde{D}(t,x)$ whenever $\tau\geq \hat{t}_1$.
Consider the event $$
\widetilde{B}_0:=\{|u_g(\hat{t}_1,x)|\leq\frac{\epsilon}{6},x\in[0,1],\quad |u_g(t,x)|\leq \frac{2\epsilon}{3}\forall t\in \hat{I}_0,x\in [0,1]\}.
$$
Then we have the following sequence of set inclusions
\begin{equation}\label{**12}
\begin{aligned}
    P(\hat{A}_0)&\geq P\left(\widetilde{B}_0\cap \{\sup_{0\leq t\leq \hat{t}_1,x\in[0,1]}|D(t,x)|\leq\frac{\epsilon}{6}\}\right)\\ &=P\left(\widetilde{B}_0\cap\{\sup_{0\leq t\leq \hat{t}_1,x\in[0,1]}|\widetilde{D}(t,x)|\leq\frac{\epsilon}{6}\}\right)\\&\geq P(\widetilde{B}_0)-P(\sup_{0\leq t\leq \hat{t}_1,x\in[0,1]}|\widetilde{D}(t,x)|>\frac{\epsilon}{6}).
\end{aligned}
\end{equation}
The equality in the second line needs some explanation. If $\tau\geq \hat{t}_1$, then $D=\widetilde{D}$ on $[0,\hat{t}_1]$. On $\widetilde{B}_0\cap \{\tau<\hat{t}_1\}$, one must have $\sup_x |u_g(\tau,x)-u_0(x)|>\epsilon$, so that $\sup_x |\widetilde{D}(\tau,x)|>\epsilon.$

Since $\beta>1$ and $\epsilon<1$, one must have $\hat{t}_1<t_1$, so that by Proposition \ref{prop2.2},
$$P(\widetilde{B}_0)\geq P(\bar{B}_0)\geq C_1\exp(-\frac{C_2}{\epsilon^2}),$$ 
where $C_1,C_2$ depend only on $\mathcal{C}_1,\mathcal{C}_2$, and where we define 
$$
\bar{B}_0:=\{|u_g({t}_1,x)|\leq\frac{\epsilon}{6},x\in[0,1],\quad |u_g(t,x)|\leq \frac{2\epsilon}{3}\forall t\in I_0,x\in [0,1]\}.
$$

It remains to estimate the probability of the last event in \eqref{**12}.
$$\begin{aligned}P(\sup_{0\leq t\leq \hat{t}_1,x\in[0,1]}|\widetilde{D}(t,x)|>\frac{\epsilon}{6})\leq\frac{1}{\sqrt{c_0}\epsilon^2}P\left(\sup_{0\leq t\leq\hat{t}_1,x\in[0,\sqrt{c_0}\epsilon^2]}|\widetilde{D}(t,x)|>\frac{\epsilon}{6}
\right)\\\leq \frac{1}{\sqrt{c_0}\epsilon^2}\frac{1}{\epsilon^{2\beta-2}}K_1\exp(-\frac{K_2}{\epsilon^{2\alpha}\mathcal{D}^2\mathcal{C}_2^2\sqrt{c_0}\epsilon^{2\beta-2}}).\end{aligned}$$
To be a bit more precise about the exponent of $\epsilon$, we take $\alpha=\epsilon^{4\beta-4}$ and the constant $\mathcal{C}:=\mathcal{D}(2\epsilon)^\alpha$ in the setting of Proposition \ref{prop2.1}.
Comparing the exponents, we see that as long as we take $\alpha+\beta>2$, we can find some $C_8,C_9$ depending only on $\mathcal{C}_1,\mathcal{C}_2$ and $\mathcal{D}$ such that 
\begin{equation}
\label{intermediate}P(\sup_{0\leq t\leq \hat{t}_1,x\in[0,1]}|\widetilde{D}(t,x)|>\frac{\epsilon}{6})\leq C_8\exp(-\frac{C_9}{\epsilon^{2(\alpha+\beta-1)}})\end{equation}
and finally 
$$\mathbb{P}(\hat{A}_0)\geq C_1\exp(-\frac{C_2}{\epsilon^2})-C_8\exp(-\frac{C_9}{\epsilon^{2(\alpha+\beta-1)}}),$$ with the first term dominating. So we conclude that when $\epsilon>0$ is sufficiently small, we can find $C_3,C_4$ depending only on $\mathcal{C}_1,\mathcal{C}_2,\mathcal{D}$ such that 
\begin{equation}\label{311113}\mathbb{P}(\hat{A_0})\geq C_3\exp(-\frac{C_4}{\epsilon^2}).\end{equation}
By the Markov property, for each $n\leq \hat{n}_1:=\lfloor\frac{T}{\epsilon^{4\beta}}\rfloor$ we have 
$$\mathbb{P}(\hat{A}_n\mid \cap_{j=0}^{n-1}\hat{A}_j)\geq C_3\exp(-\frac{C_4}{\epsilon^2}).$$
Thus for some constant $C_1>0$ depending on $\mathcal{C}_1,\mathcal{C}_2$,$$P\left(|u(t,x)|\leq\epsilon,t\in[0,T],x\in[0,1]\right)\geq \exp(-\frac{C_1}{\epsilon^{2}}\frac{T}{\epsilon^{4\beta}})\geq \exp(-\frac{C_1T}{\epsilon^{4\beta+2}}).$$
The various constants $C_1,C_2,\cdots,C_8,C_9$ depend only on $\mathcal{C}_1,\mathcal{C}_2$ and may change from line to line. This establishes the lower bound in \eqref{mainest}.

Finally, we note that we may take a different estimate of \eqref{intermediate} and deduce that, for some constants $C_8',C_9'$ depending on $\mathcal{C}_1,$ $\mathcal{C}_2$ but not $\mathcal{D}$, 
$$P(\sup_{0\leq t\leq \hat{t}_1,x\in[0,1]}|\widetilde{D}(t,x)|>\frac{\epsilon}{6})\leq C_8'\exp(-\frac{C_9'}{\epsilon^{2(\alpha+\beta-1)}\mathcal{D}^2}).$$
Then there exists some $\mathcal{D}_0$ depending only on $\mathcal{C}_1$, $\mathcal{C}_2$ so that, for any $\mathcal{D}<\mathcal{D}_0$ and $\beta=2-\alpha$, we can find $C_3'$, $C_4'$ depending on $\mathcal{C}_1$, $\mathcal{C}_2$ such that (note that \eqref{311113} requires $\beta+\alpha>2$)
$$\mathbb{P}(\hat{A_0})\geq C_3'\exp(-\frac{C_4'}{\epsilon^2}).$$
The rest of the estimate proceeds exactly the same as before. The net benefit we get from this business is that whenever $\mathcal{D}<\mathcal{D}_0$, all the estimates in Theorem \ref{theorem1.1} can be done in the case $\alpha+\beta=2$. This justifies the last assertion of Theorem \ref{theorem1.1}.

\section*{Statements and Declarations}
The author receives financial support from EPSRC grant EP/W524141/1.

The author has no competing interests to declare that are relevant to the content of this article.

Data sharing not applicable to this article as no datasets were generated or analysed during the current study.

\printbibliography
\end{document}